\documentclass[12pt, reqno]{amsart}

\usepackage{amsmath,amsthm,amssymb}
\usepackage{hyperref}

\newtheorem{theorem}{Theorem}
\newtheorem{corollary}[theorem]{Corollary}
\newtheorem{lemma}[theorem]{Lemma}

\newcommand{\C}{\mathbb{C}}
\newcommand{\R}{\mathbb{R}}
\newcommand{\Z}{\mathbb{Z}}
\newcommand{\N}{\mathbb{N}}
\newcommand{\T}{{\mathbb T}}

\newcommand{\supp}{\operatorname{supp}}
\newcommand{\re}{\operatorname{Re}}
\newcommand{\im}{\operatorname{Im}}

\newcommand{\e}{\varepsilon}

\newcommand{\cH}{{\mathcal H}}

\newcommand{\be}{\begin{equation}}
\newcommand{\ee}{\end{equation}}

\begin{document}

\title[Local $L^\infty$ Bounds via Diophantine Problems]{Local $L^\infty$ Bounds for Eigenfunctions \\ of Complex Elliptic Operators via Diophantine Problems}

\author{Omer Friedland} 
\address{Institut de Math\'ematiques de Jussieu, Sorbonne Universit\'e, 4 place Jussieu, 75005 Paris}
\email{omer.friedland@imj-prg.fr}

\author{Henrik Uebersch\"ar}
\address{Institut de Math\'ematiques de Jussieu, Sorbonne Universit\'e, 4 place Jussieu, 75005 Paris}
\email{henrik.ueberschar@imj-prg.fr}



\thanks{We thank Y. Yomdin and N. Lerner for valuable comments on an earlier draft. We are especially grateful to N. Lerner for insightful discussions on constant-coefficient hypoelliptic operators and for suggestions that improved the exposition.}

\begin{abstract}
We prove local bounds on the amplitude of  eigenfunctions of complex constant-coefficient elliptic operators with a smooth potential on an arbitrary open subset of $\R^d$ by estimating it in terms of the number of solutions of a diophantine inequality arising from the symbol of the operator. 
In the special case of positive elliptic operators, we recover H\"ormander's classical exponent up to an arbitrarily small loss. We show that a much better exponent may be obtained when the principal symbol of the operator has complex coefficients. We generalize our estimate to any higher-order derivatives of eigenfunctions.


\end{abstract}

\maketitle

\section{Introduction}

In this article, we establish a pointwise bound on the amplitude of the eigenfunctions of an elliptic operator $H+V$ on an open set $\Omega\subset\R^d$, where $H$ has constant complex coefficients and $V\in C^\infty(\Omega,\C)$. We introduce a new method which reduces the problem to estimating the lattice point count in shells around the $\lambda$-level set of the principal symbol of the operator. While such estimates are classical in the case of positive elliptic operators and intimately tied to the local Weyl law and the spectral function's fine behavior near the diagonal \cite{Le1952,Av1956,Hormander1968,Agmon1968,DG1975,Seeley1978}, very little is known about the general case of complex coefficients.



Our main contributions may be summarized as follows.
\begin{enumerate}
\item We prove local bounds on the amplitude of solutions of
$$(H+V)\psi_\lambda = \lambda\psi_\lambda, \quad \lambda\in\C,$$ as well as bounds on their derivatives, on arbitrary open sets
$\Omega\subset\mathbb R^d$, where $H$ is a complex constant-coefficient
elliptic operator and $V\in C^\infty(\Omega,\mathbb C)$.\\

\item These bounds are expressed purely in terms of a Diophantine
lattice-point counting function $\mathcal N_\lambda$ in frequency space,
measuring the thickness of a shell around the level set $\{P(\xi)=\lambda\}$.\\

\item In the classical case of positive elliptic operators with real coefficients,
our method recovers H\"ormander’s sharp exponent $(d-1)/(2m)$ up to an arbitrarily
small loss (see Subsection~\ref{comparison}). For broad families of complex symbols, however,
we show that $\mathcal N_\lambda$ is much smaller and obtain strictly better exponents
(down to $(d-2)/(2m)+\varepsilon$ in balanced regimes).\\
\end{enumerate}

A natural question to ask is, how our general bound ({\em valid for any complex coefficient operator}) compares with the classical bounds which are available in the {\em  special case of positive elliptic operators}. We show that the exponent of the eigenvalue which appears in our bound depends on the type of complex elliptic operator that is chosen. We illustrate this in the case of $2$-dimensional domains, where we show that for certain operators our bound gives a much better exponent than H\"ormander's classical exponent (see subsection \ref{comparison} for a detailed discussion).

\subsection{Main Results}

We consider smooth solutions \footnote{Note that this holds in particular if $H+V-\lambda$ is hypoelliptic, notably for elliptic operators with smooth coefficients.} of the following eigenvalue problem on an open domain $\Omega \subset \R^d$
\be \label{eq:main}
(H + V)\psi_\lambda = \lambda \psi_\lambda \quad \text{on } \Omega,
\ee
where $V\in C^\infty(\Omega,\C)$ and $\lambda\in\C$ with $|\lambda|\geq1$.

Let $H$ be an $m$th-order partial differential operator with constant coefficients
$$
H = \sum_{|\alpha|\leq m}c_\alpha D^\alpha, \quad c_\alpha\in\C,
$$
where we denote by $\alpha = (\alpha_1,\cdots,\alpha_d)\in\N^d$ a multi-index and use the standard notation $D^\alpha = (-i\partial_{x_1})^{\alpha_1}\cdots (-i\partial_{x_d})^{\alpha_d}$. The total order of the multi-index is denoted by $|\alpha| = \alpha_1+\cdots+\alpha_d$.

We denote the symbol of $H$ by 
$$
P(\xi) = \sum_{|\alpha|\leq m}c_\alpha \xi^\alpha, 
$$
where we use the convention $\xi^\alpha = \xi_1^{\alpha_1} \cdots \xi_d^{\alpha_d}$ for any $\xi\in\R^d$. 

Moreover, we denote the principal symbol of $H$ by
$$
P_m(\xi) = \sum_{|\alpha|=m}c_\alpha \xi^\alpha.
$$

We denote by $\mathcal{N}_\lambda$ the number of lattice points in a shell around the level set $\{P(\xi) = \lambda\}$ which is defined by a Diophantine inequality:
\[
\mathcal{N}_\lambda : = \#\left\{ \xi \in \mathbb{Z}^d : |P(\xi) - \lambda| \le |\xi|^{m-1+\delta} \right\}.
\]
Although it depends on the symbol $P(\xi)$ and on $\delta>0$, we suppress this dependence in the notation for simplicity.

Our main estimate is the following local bound on derivatives of eigenfunction which we prove in section \ref{sec:local-bound}.

\begin{theorem}[Local Bound] \label{thm:local-bound}
Let \( \delta,r > 0 \) and $\gamma$ be any multi-index. Then for any ball \( B(x,r) \subset \Omega \), we have
\begin{equation}\label{local-gen}
|D^\gamma\psi_\lambda(x)| \lesssim_{d,r,\delta} |\lambda|^{|\gamma|/m}\left(1 + \mathcal{N}_\lambda^{1/2} \right) \|\psi_\lambda\|_{L^2(B(x,r))}.
\end{equation}
Here, the implied constant depends on $d,r,\delta$.
\end{theorem}

We stress that our estimate is still valid if one allows the radius $r$ to shrink slowly in the sense that $r^{-1}\lesssim(\re\lambda)^{o(1)}$, as $\re\lambda\to+\infty$. We make this explicit only in the case of Schr\"odinger operators, however, this refinement may be justified for any constant coefficient operator $H$ (for more details see Section~\ref{sec:shrink}).

We obtain the following global interior bound as a direct corollary of the local estimate \eqref{local-gen}. This recovers the sharp exponent from H\"ormander’s classical estimate, up to $\varepsilon$, in a vastly more general setting: it applies to arbitrary open domains, including those with rough or fractal boundaries, and requires no assumptions on compactness, boundary conditions, or self-adjointness (for example, $V$ can be complex-valued).
\begin{corollary}[$\varepsilon$-sharp Interior $L^\infty$ Bound] \label{cor:interior-bound}
Let $H$ be an elliptic operator of order $m \geq 1$ with constant coefficients, $V \in C^\infty(\Omega, \mathbb{C})$, and $\psi_\lambda$ a solution to $(H + V)\psi_\lambda = \lambda \psi_\lambda$ on an open set $\Omega \subset \mathbb{R}^d$. For $r>0$, we define $\Omega_{-r}:=\{x\in\Omega \mid B(x,r)\subset\Omega\}$.

Then for every $\varepsilon > 0$ and any $r > 0$,
\[
\|D^\gamma\psi_\lambda\|_{L^\infty(\Omega_{-r})} \lesssim_{\varepsilon, r} |\lambda|^{|\gamma|/m}\left(1 + \mathcal{N}_\lambda^{1/2} \right) \|\psi_\lambda\|_{L^2(\Omega)},
\]
uniformly as $|\lambda| \to \infty$. The implicit constant depends on $H$, $d$, $r$, $\varepsilon$, and bounds on derivatives of $V$.
\end{corollary}

\subsection{Comparison to Existing Results}\label{comparison}

First, we compare our main estimate with the known bound in the case of positive elliptic operators to affirm that our bound is $\e$-sharp. Then we show how one may beat this exponent in the case of complex coefficients.\\

{\bf $\varepsilon$-sharp exponent for real coefficients.}
Suppose \( H \) is elliptic of order \( m \) with {\em real} coefficients. Then, for any \( \delta > 0 \) (cf. Lemma \ref{lem:lattice-count} and its proof in section \ref{sec:count}), we have the following estimate on the lattice point count:
$$\mathcal{N}_\lambda\lesssim_{d,r,\delta}|\lambda|^{(d-1+\delta)/m}.$$
In combination with Theorem \ref{thm:local-bound} this yields the following local bound which recovers H\"ormander's exponent up to $\varepsilon$:
\begin{equation}\label{local-real}
|\psi_\lambda(x)| \lesssim_{d,r,\e} |\lambda|^{(d-1)/(2m)+\varepsilon} \|\psi_\lambda\|_{L^2(B(x,r))}.
\end{equation}
In particular, this shows that, up to the $\varepsilon$-loss in the exponent, the bound in Theorem \ref{thm:local-bound} is optimal.\\

{\bf Better exponent for complex coefficients.} 
We stress that one may obtain much better bounds if the principal symbol of the operator $H$ has complex coefficients. For example, choose $H=-\partial_{x}^2-i\partial_{y}^2$. In this case, we have  the following bound which is a consequence of a better estimate on $\mathcal{N}_\lambda$ in this case (cf. Lemma \ref{lem:complex-count} and its proof in section \ref{sec:count})
\begin{equation}
|\psi_\lambda(x)|\lesssim_\e F_\e(\lambda)\|\psi_\lambda\|_{L^2(B(x,r)}
\end{equation}
where 
\begin{equation}
\begin{split}
F_\e(\lambda)=
\begin{cases} 
(\re\lambda)^{1/4+\e}(\im\lambda)^{-1/4}\quad &\text{if}\; (\re\lambda)^{1/2}\leq\im\lambda\leq\re\lambda,\\
\\
(\re\lambda)^{1/8+\e} \quad &\text{if}\; \im\lambda\leq(\re\lambda)^{1/2},
\end{cases}
\end{split}
\end{equation}
if $0<\im\lambda\leq\re\lambda\to+\infty$. Analogous bounds hold if one exchanges the roles of $\re\lambda$ and $\im\lambda$. If we recall that H\"ormander's bound for the Laplacian on $2$-dimensional domains gives $O(\lambda^{1/4}$), we clearly observe that for the general case of complex coefficients much better bounds can be obtained.

We note that it is easy to find a general class of examples that give an improvement over H\"ormander's exponent in the complex case. For example one might take $P(\xi)=Q_1(\xi_1,\cdots,\xi_p)+iQ_2(\xi_{p+1},\cdots,\xi_d)$, where $$Q_1(\xi_1,\cdots,\xi_p)=\sum_{j=1}^p \alpha_j \xi_j^2,$$ and $$Q_2(\xi_{p+1},\cdots,\xi_d)=\sum_{j=p+1}^d \alpha_j \xi_j^2,$$ for $\alpha_j>0$, $j=1,\cdots,d$. In the regime $\re\lambda\asymp \im\lambda\to+\infty$, this gives an upper bound of the form $F_{\e,d}(\lambda)=|\lambda|^{(d-2)/(2m)+\e}$ which still beats H\"ormander's exponent $(d-1)/(2m)$.

Moreover, we could even take $Q_1$, $Q_2$ which give rise to hyperbolic lattice point problems: $Q_1(\xi_1,\xi_2)=\xi_1\xi_2$, $Q_2(\xi_1,\xi_2)=\xi_1^2-\xi_2^2$, where $P=Q_1+iQ_2$ gives rise to an elliptic operator. 



\subsection{Concluding remarks}

Our approach is based on taking a smooth cutoff of the original eigenfunction and imbedding this approximate eigenfunction in a fixed $d$-torus. We then estimate the Fourier coefficients of these approximate eigenfunctions. The summation in the frequency domain is truncated outside a sufficiently thick neighborhood of the level set \( \{ \xi \in \mathbb{R}^d : P(\xi) = \lambda \} \), where \( P(\xi) \) is the symbol of the operator. The summation over the lattice points in this neighbourhood then gives the main term in our estimate. Moreover, our approach also works for derivatives of eigenfunctions and can be adapted to shrinking balls, as long as $r^{-1}=|\lambda|^{o(1)}$.

Our method allows us to obtain local estimates of eigenfunctions and their derivatives for elliptic operators whose principal symbol has complex coefficients. This goes far beyond classical approaches for positive elliptic operators which rely on the validity of a spectral theorem and the connection of such estimates with the local Weyl law.

We expect our method to offer a new bridge between Spectral Geometry and Diophantine Geometry  --  particularly in the case of operators and domains where classical methods fail as in the case of elliptic operators with complex coefficients and domains with fractal boundary. In forthcoming work, we will apply the method to the study of nodal sets for such operators.

\subsection{Structure of the paper}

The remainder of the paper is organized as follows:

In Section~\ref{sec:local-bound}, we prove the main local $L^\infty$ estimate (Theorem~\ref{thm:local-bound}) by reducing pointwise control of eigenfunctions to a Diophantine lattice counting problem. The argument relies on a local Fourier decomposition on the torus and precise bounds on Fourier coefficients of localized eigenfunctions.

Section~\ref{sec:count} is devoted to analyzing the lattice point counting problem. We establish the connection between ellipticity and the finiteness of Diophantine shells, prove asymptotic estimates for the number of integer solutions, and show how homogeneity of the symbol controls the geometry of frequency shells.


Appendix~\ref{sec:rdn} reviews classical estimates on the representation of integers as sums of squares, used in bounding sums over lattice points. Appendix~\ref{sec:shrink} refines our estimates to include the dependence of constants on the shrinking radius $r$, and shows that the $L^\infty$ bounds remain $\e$-sharp even when $r = r(\lambda)$ decays slowly, as $|\lambda|\to+\infty$.

\section{Proof of Theorem ~\ref{thm:local-bound}} \label{sec:local-bound}

We fix $\chi\in C^\infty_c(\R^d)$ such that $\supp\chi = B(0,1)$, it satisfies $0\le \chi\le 1$ and $\chi(0) = 1$. Given $r > 0$ and $x_0\in \Omega_{-r}$, we define $\chi_r(x) = \chi\left(\frac{x-x_0}{r} \right)$. Introducing a smooth cutoff of $\psi_\lambda$ as $\widetilde{\psi}_\lambda = \chi_r \psi_\lambda$, we proceed to define a function 
$$
\Psi(x) = \sum_{k\in 2\pi\Z^d} \widetilde{\psi}_\lambda(x-k)
$$
on $\T^d = \R^d/2\pi\Z^d$. Let us denote by $\cH^d$ the hypercube with vertices $x_0+\pi(\pm1,\cdots,\pm1)$ which we take as a fundamental domain. 
In particular, we have for any $x\in\cH^d$ that 
$$
\Psi(x) = \chi_r(x)\psi_\lambda(x)
$$
because the support of the translates $\tilde\psi_\lambda(\cdot-k)$, for $k\in 2\pi\Z^d$, namely the balls $B(x_0+k,r)$, are disjoint due to the assumption $r<1$.

We expand $\Psi$ in its Fourier series on $\T^d$ 
$$
\Psi(x) = \sum_{\xi\in\Z^d} \hat{\Psi}(\xi)e_\xi(x), \quad e_\xi(x): = e^{i\xi\cdot x}.
$$
We have for any $x\in\T^d$,
\be \label{eq:sum}
|\Psi(x)| \le \sum_{\xi\in\Z^d} |\hat{\Psi}(\xi)|
\ee
where we note that the Fourier coefficients decay rapidly, because $\Psi$ is smooth.

We have the following bound on the Fourier coefficients of $\Psi$ (refer to Section ~\ref{sec:fourier-bound} for its proof) which we will use in the proof:

\begin{lemma} \label{lem:fourier-bound}
Let $\alpha\in\N$. For any $\xi\in\Z^d\setminus\{0\}$ such that $P(\xi)\neq\lambda$, we have
\be \label{eq:fourier-bound}
|\hat{\Psi}(\xi)| \lesssim_{\alpha,r} \frac{|\xi|^{(m-1)\alpha}}{|P(\xi)-\lambda|^{\alpha}} \|\psi_\lambda\|_{L^2(B(x_0,r))}.
\ee
\end{lemma}
Our objective is to derive a bound on the sum \eqref{eq:sum}. 
Let us fix $\delta\in(0,1)$ and take $\alpha = \alpha(\delta)>\frac{d}{\delta}$. We proceed by splitting the summation $\sum_{\xi \in \mathbb{Z}^d} |\hat{\Psi}(\xi)|$ into two parts.

\subsection{Part I: $|P(\xi)-\lambda| < |\xi|^{m-1+\delta}$}

We use the Cauchy-Schwarz inequality to obtain
\begin{align*}
\sum_{|P(\xi)-\lambda| < |\xi|^{m-1+\delta}} |\hat{\Psi}(\xi)| \le & (\sum_{|P(\xi)-\lambda| < |\xi|^{m-1+\delta}} |\hat{\Psi}(\xi)|^2)^{1/2}(\sum_{|P(\xi)-\lambda| < |\xi|^{m-1+\delta}}1)^{1/2} \\
\leq& \|\Psi\|_2 (\sum_{|P(\xi)-\lambda| < |\xi|^{m-1+\delta}}1)^{1/2} 
\end{align*}

Thus, we obtain the following bound on part I of the summation:
\begin{align*}
\sum_{|P(\xi)-\lambda| < |\xi|^{m-1+\delta}} |\hat{\Psi}(\xi)| 
\lesssim_{r,d,\delta} \|\psi_\lambda\|_{L^2(B(x_0,r))} {\mathcal N}_\lambda^{1/2}
\end{align*}
where we used $\|\Psi\|_2\leq\|\psi\|_{L^2(B(x_0,r))}$.

\subsection{Part II: $|P(\xi)-\lambda| \ge|\xi|^{m-1+\delta}$}

First of all, we have
$$
|\hat{\Psi}(0)|\lesssim \|\psi_\lambda\|_{L^2(B(x_0,r))}.
$$

For $\xi\neq0$, we recall the bound from Lemma ~\ref{lem:fourier-bound}
$$
|\hat{\Psi}(\xi)|\lesssim \frac{|\xi|^{(m-1)\alpha}}{|P(\xi)-\lambda|^\alpha} \|\psi_\lambda\|_{L^2(B(x_0,r))}.
$$

Using this bound on the Fourier coefficients, we obtain
\begin{align*}
 \sum_{|P(\xi)-\lambda| \ge |\xi|^{m-1+\delta}} |\hat{\Psi}(\xi)| &\lesssim_{} \|\psi_\lambda\|_{L^2(B(x_0,r))} \sum_{|P(\xi)-\lambda| \ge |\xi|^{m-1+\delta}} \frac{|\xi|^{(m-1)\alpha}} {|P(\xi)-\lambda|^\alpha} \\
& \le \|\psi_\lambda\|_{L^2(B(x_0,r))} \sum_{\xi\in\Z^d\setminus\{0\}} |\xi|^{-\alpha\delta}.
\end{align*}

Recall $\delta > d/\alpha$. Then, we get (set $\e'' = \alpha\delta-d$)
\begin{align*}
\sum_{|P(\xi)-\lambda|\ge|\xi|^{m-1+\delta}} |\hat{\Psi}(\xi)| 
& \lesssim \|\psi_\lambda\|_{L^2(B(x_0,r))} \sum_{n \ge 1} r_d(n)n^{-(d+\e'')/2} \\
& = \|\psi_\lambda\|_{L^2(B(x_0,r))} \zeta_d\left(\frac{d+\e''}{2} \right), 
\end{align*}
where 
$$
r_d(n): = \#\{\xi\in\Z^d \mid n = |\xi|^2\}
$$
and 
$$
\zeta_d(s) = \sum_{n \ge 1}r_d(n)n^{-s},
$$
which converges absolutely in the half-plane $\re s>d/2$ in view of the bound 
$r_d(n)\lesssim n^{d/2-1+o(1)}$ (see Appendix ~\ref{sec:rdn} for a review of estimates on the function $r_d(n)$).

\subsection{Proof of Lemma ~\ref{lem:fourier-bound}} \label{sec:fourier-bound}

First we need the following simple observation:

\begin{lemma} \label{lem-deriv-bound}
Let $\gamma = (\gamma_1,\cdots,\gamma_d) \in\N^d$ be a multi-index and denote $D^\gamma = (-i\partial_{x_1})^{\gamma_1} \cdots(-i\partial_{x_d})^{\gamma_d}$.
Let $h\in C_c^\infty$ be a cutoff that imbeds into $\T^d$. We have the bound
$$
\left|\int_{\T^d}h \cdot (D^\gamma\psi_\lambda) e_{-\xi}d\mu\right| \lesssim_{\gamma,h,d} |\xi|^{|\gamma|} \|\psi_\lambda\|_{L^2(\supp h)}.
$$
\end{lemma}

\begin{proof}
An integration by parts yields
\begin{align*}
& \left|\int_{\T^d}h \cdot (D^\gamma\psi_\lambda) e_{-\xi}d\mu\right| = \left|\int_{\T^d} \psi_\lambda D^\gamma(he_{-\xi})d\mu\right| \\
& = \left|\int_{\supp h} \psi_\lambda D^\gamma(he_{-\xi})d\mu\right|\le \|\psi_\lambda\|_{L^2(\supp h)} \cdot \|D^\gamma(he_{-\xi})\|_2.
\end{align*}

The derivative $D^\gamma(he_{-\xi})$ consists of terms which are products of derivatives of $h$ and $e_{-\xi}$, where the sum of the total order equals $|\gamma|$,
$$
D^\gamma(he_{-\xi}) = \sum_{\gamma_1+\gamma_2 = \gamma}C(\gamma_1,\gamma_2)(D^{\gamma_1}h)\cdot(D^{\gamma_2}e_{-\xi}),
$$
and the coefficients $C(\gamma_1,\gamma_2)$ are products of binomial coefficients.

Moreover, 
$$
D^{\gamma'}(e_{-\xi}) = (-1)^{|\gamma'|} \xi^{\gamma'}e_{-\xi}
$$
and
$$
\|D^{\gamma'}(e_{-\xi})\|_2\le (2\pi)^{d/2} |\xi|^{|\gamma'|}.
$$

So, for $|\xi|\geq1$, we have 
\begin{align*}
\left|\int_{\T^d}h (D^\gamma\psi_\lambda) e_{-\xi}d\mu\right| 
&\leq \|\psi_\lambda\|_{L^2(\supp h)} \sum_{\gamma_1+\gamma_2 = \gamma}|C(\gamma_1,\gamma_2)|\|D^{\gamma_1}h\|_2\|D^{\gamma_2}e_{-\xi}\|_2\\
&\lesssim_{h,\delta,|\gamma|} |\xi|^{|\gamma|} \|\psi_\lambda\|_{L^2(\supp h)}.
\end{align*}
\end{proof}

Let 
$$
H = \sum_{|\gamma|\leq m} c_\gamma D^\gamma
$$
be an $m$th order partial differential operator on $\Omega$. Our claim is that when we expand the term $(H-\lambda)^\alpha\Psi$ into derivatives of $\psi_\lambda$, the highest order derivative that will appear is of order $\leq \alpha(m-1)$. We prove this by induction on $\alpha$.
\\

{\bf The case $\alpha = 1$.}

First of all, recall 
\begin{align*}
D^\gamma(\chi \psi_\lambda) = &\sum_{\gamma_1+\gamma_2 = \gamma}C(\gamma_1,\gamma_2)D^{\gamma_1}\chi D^{\gamma_2}\psi_\lambda\\
 = &\sum_{\substack{\gamma_1+\gamma_2 = \gamma \\ \gamma_1\neq 0}} C(\gamma_1,\gamma_2)D^{\gamma_1}\chi D^{\gamma_2}\psi_\lambda+\chi D^\gamma \psi_\lambda
\end{align*}
where we recall that the coefficients $C(\gamma_1,\gamma_2)$ are products of binomial coefficients and $C(0,\gamma) = 1$.

We may now use this identity to compute
\begin{align*}
(H-\lambda) \Psi & = H(\chi\psi_\lambda)-\lambda\chi\psi_\lambda \\
& = \sum_{|\gamma|\leq m}c_\gamma \sum_{\substack{\gamma_1+\gamma_2 = \gamma \\ \gamma_1\neq0}}C(\gamma_1,\gamma_2)(D^{\gamma_1}\chi)(D^{\gamma_2}\psi_\lambda)+\chi(H-\lambda) \psi_\lambda \\
& = \sum_{|\gamma|\leq m}c_\gamma \sum_{\substack{\gamma_1+\gamma_2 = \gamma \\ \gamma_1\neq0}}C(\gamma_1,\gamma_2)(D^{\gamma_1}\chi)(D^{\gamma_2}\psi_\lambda)
-\chi V\psi_\lambda.
\end{align*}
And the highest order derivative on $\psi_\lambda$ is of order $\leq m-1$, because $\gamma_1\neq 0$ and thus $|\gamma_2|<|\gamma|\leq m$.

Thus, the claim holds true when $\alpha = 1$.
\\

{\bf Induction step.}

Assuming the claim holds for $\alpha\ge1$, we demonstrate it for step $\alpha+1$, as follows 
$$
(H-\lambda)^{\alpha+1} \Psi = (H-\lambda)(H-\lambda)^{\alpha} \Psi.
$$

By the induction hypothesis, we have that the highest order derivative of $\psi_\lambda$ which appears in the interior term is of order $\leq \alpha(m-1)$. The highest order term is of the form $hD^\eta\psi_\lambda$, where $\eta$ is a multi-index with $|\eta| \leq \alpha(m-1)$ and $h$ is a smooth cutoff. 

Let us calculate 
\begin{align*}
& (H-\lambda)(hD^\eta\psi_\lambda) = \\
 = & \sum_{\substack{|\gamma|\leq m  \\ \gamma_1+\gamma_2 = \gamma \\ \gamma_1\neq0}}c_\gamma C(\gamma_1,\gamma_2)(D^{\gamma_1}h)D^{\gamma_2}(D^\eta\psi_\lambda)+h(H-\lambda) D^\eta\psi_\lambda \\
 = & \sum_{\substack{|\gamma|\leq m \\  \gamma_1+\gamma_2 = \gamma \\ \gamma_1\neq0}}c_\gamma C(\gamma_1,\gamma_2)(D^{\gamma_1}h)D^{\gamma_2}(D^\eta\psi_\lambda)+hD^\eta(H-\lambda)\psi_\lambda 
\end{align*}
which we rewrite as follows, using our assumption that $\psi_\lambda$ is a solution of \eqref{eq:main},
\begin{align*}
& \sum_{\substack{|\gamma|\leq m \\  \gamma_1+\gamma_2 = \gamma \\ \gamma_1\neq0}}c_\gamma C(\gamma_1,\gamma_2)(D^{\gamma_1}h)(D^{\gamma_2+\eta}\psi_\lambda)-hD^\eta(V\psi_\lambda)\\
 = & \sum_{\substack{|\gamma|\leq m \\  \gamma_1+\gamma_2 = \gamma  \\ \gamma_1\neq0}}c_\gamma C(\gamma_1,\gamma_2)(D^{\gamma_1}h)(D^{\gamma_2+\eta}\psi_\lambda)\\
&\quad -\sum_{\eta_1+\eta_2 = \eta}C(\eta_1,\eta_2)h(D^{\eta_1}V)(D^{\eta_2}\psi_\lambda)
\end{align*}
In the sum, the highest order derivative on $\psi_\lambda$ is of order $\leq (m-1)+|\eta|$, because $\gamma_1\neq0$ which means $|\gamma_2|<|\gamma|\leq m$.
Moreover, $|\eta| \leq \alpha(m-1)$, which yields that the highest order derivative is of order $\leq (\alpha+1)(m-1)$. This concludes the induction step.
\\

{\bf Bounding the Fourier coefficient.}

To finish the proof we apply the Fourier transform
$$
(P(\xi)-\lambda)^\alpha\hat{\Psi}(\xi) = \int_{\T^d}[(H-\lambda)^\alpha\Psi]e_{-\xi}d\mu
$$
and note that Lemma ~\ref{lem-deriv-bound} yields
$$
\left|\int_{\T^d}[(H-\lambda)^\alpha\Psi]e_{-\xi}d\mu\right| \lesssim_{r,\alpha} |\xi|^{(m-1)\alpha} \cdot \|\psi_\lambda\|_{L^2(B(x_0,r))},
$$
because $(H-\lambda)^\alpha\Psi$ can be written as a sum of products of derivatives of $\chi_r$ and $V$ and derivatives of $\psi_\lambda$, the former of which are supported on $B(x_0,r)$ and the latter of which are of order at most $(m-1)\alpha$. Combining all of the above, we obtain 
\begin{align*}
|P(\xi)-\lambda|^\alpha |\hat{\Psi}(\xi)| \lesssim_{r,\alpha} |\xi|^{(m-1)\alpha} \|\psi_\lambda\|_{L^2(B(x_0,r))},
\end{align*}
which concludes the proof of Lemma ~\ref{lem:fourier-bound}.

\subsection{Generalization to derivatives}\label{sec:deriv}

Let $\chi \in C^\infty_c(\mathbb{R}^d)$ be such that $\chi \equiv 1$ on $B(0,1/2)$ and $\operatorname{supp} \chi \subset B(0,1)$. For $r > 0$ and $x_0 \in \Omega_{-r}$, define the localized cutoff $\chi_r(x) := \chi\left( \frac{x - x_0}{r} \right)$. We define the periodized function
\[
\Psi(x) = \sum_{k \in 2\pi \mathbb{Z}^d} (\chi_r \psi_\lambda)(x - k),
\]
as in the proof of Theorem~\ref{thm:local-bound}. Then $\Psi = \psi_\lambda$ in $B(x_0, r/2)$, and the same holds for all derivatives: $D^\gamma \Psi = D^\gamma \psi_\lambda$ on $B(x_0, r/2)$.

Using the Fourier identity $\widehat{D^\gamma \Psi}(\xi) = \xi^\gamma \hat{\Psi}(\xi)$, we estimate
\[
|D^\gamma \psi_\lambda(x_0)| = |D^\gamma \Psi(x_0)| \le \sum_{\xi \in \mathbb{Z}^d} |\xi|^{|\gamma|} |\hat{\Psi}(\xi)|.
\]

We now split the summation into two parts as before.

{\em Part I: The Diophantine Shell}

Consider the set
\[
S_\lambda := \left\{ \xi \in \mathbb{Z}^d : |P(\xi) - \lambda| \le |\xi|^{m - 1 + \delta} \right\}.
\]
For $\xi \in S_\lambda$, we estimate
\[
|\xi|^{|\gamma|} |\hat{\Psi}(\xi)| \le |\xi|^{|\gamma|} \cdot \|\Psi\|_2 \lesssim \lambda^{\frac{|\gamma|}{m}} \|\psi_\lambda\|_{L^2(B(x_0, r))},
\]
since in the shell we have \( |\xi| \asymp |\lambda|^{1/m} \). Since $\# S_\lambda = \mathcal{N}_\lambda$, an application of Cauchy-Schwarz gives
\[
\sum_{\xi \in S_\lambda} |\xi|^{|\gamma|} |\hat{\Psi}(\xi)| \lesssim \lambda^{\frac{|\gamma|}{m}}(1+\mathcal{N}_\lambda^{1/2}) \|\psi_\lambda\|_{L^2(B(x_0, r))}.
\]

{\em Part II: Complement of the Shell}

For $\xi \notin S_\lambda$, we use Lemma~\ref{lem:fourier-bound}:
\[
|\hat{\Psi}(\xi)| \lesssim \frac{|\xi|^{(m - 1)\alpha}}{|P(\xi) - \lambda|^\alpha} \|\psi_\lambda\|_{L^2(B(x_0, r))},
\]
and hence,
\[
|\xi|^{|\gamma|} |\hat{\Psi}(\xi)| \lesssim \frac{|\xi|^{(m - 1)\alpha + |\gamma|}}{|P(\xi) - \lambda|^\alpha} \|\psi_\lambda\|_{L^2(B(x_0, r))}.
\]
Outside the shell, \( |P(\xi) - \lambda| \ge |\xi|^{m - 1 + \delta} \), so
\[
|\xi|^{|\gamma|} |\hat{\Psi}(\xi)| \lesssim |\xi|^{|\gamma| - \alpha \delta} \|\psi_\lambda\|_{L^2(B(x_0, r))}.
\]
Summing over \( \xi \in \mathbb{Z}^d \) gives convergence provided \( \alpha \delta > |\gamma| + d \), i.e.,
\[
\sum_{\xi \notin S_\lambda} |\xi|^{|\gamma|- \alpha \delta} |\hat{\Psi}(\xi)| \lesssim \|\psi_\lambda\|_{L^2(B(x_0, r))}.
\]

{\em Conclusion}

Combining both contributions, we obtain:
\[
|D^\gamma \psi_\lambda(x_0)| \lesssim \lambda^{\frac{|\gamma|}{m}}(1+\mathcal{N}_\lambda^{1/2}) \|\psi_\lambda\|_{L^2(B(x_0, r))}.
\]

\section{Lattice Point Counting and Ellipticity} \label{sec:count}

In this section, we present several lemmas that make the connection between ellipticity and our lattice point counting problem. 

\begin{lemma}[Finiteness of Diophantine Shells and Ellipticity] \label{lem:ellipticity-finiteness}
Let \( H = \sum_{|\gamma| \leq m} c_\gamma D^\gamma \) be a constant-coefficient partial differential operator on \( \mathbb{R}^d \), with complex coefficients \( c_\gamma \), and let \( P(\xi) \) denote its symbol. Fix \( \delta > 0 \) and \( \lambda \in \mathbb{C} \). Then the following are equivalent:
\begin{enumerate}
 \item \( H \) is elliptic, i.e., its principal symbol \( \sigma(\xi) = \sum_{|\gamma| = m} c_\gamma \xi^\gamma \) satisfies \( \sigma(\xi) \neq 0 \) for all \( \xi \in \mathbb{R}^d \setminus \{0\} \),
 \item The Diophantine inequality
 \[
 |P(\xi) - \lambda| \leq |\xi|^{m-1+\delta}
 \]
 admits only finitely many integer solutions \( \xi \in \mathbb{Z}^d \).
\end{enumerate}
\end{lemma}

\begin{proof}
We prove both implications.

\textbf{(1) implies (2):} 
Write the symbol as \( P(\xi) = \sigma(\xi) + r(\xi) \), where \( \sigma \) is the homogeneous principal part of degree \( m \), and \( r(\xi) \) contains lower-order terms, so that \( |r(\xi)| \lesssim |\xi|^{m-1} \) as \( |\xi| \to \infty \).

Ellipticity implies that there exists \( C > 0 \) such that \( |\sigma(\xi)| \geq C|\xi|^m \) for all \( \xi \in \mathbb{R}^d \) with \( |\xi| \) sufficiently large. Then:
\[
|P(\xi) - \lambda| \geq |\sigma(\xi)| - |r(\xi)| - |\lambda| \gtrsim |\xi|^m - C'|\xi|^{m-1},
\]
for some constant \( C' \), so that the inequality
\[
|P(\xi) - \lambda| \leq |\xi|^{m-1+\delta}
\]
can only be satisfied for finitely many \( \xi \in \mathbb{Z}^d \), completing this direction.

\textbf{(2) implies (1):} 
Assume \( H \) is not elliptic. Then there exists \( \zeta \in \mathbb{R}^d \setminus \{0\} \) such that \( \sigma(\zeta) = 0 \). By homogeneity, \( \sigma(t\zeta) = t^m \sigma(\zeta) = 0 \) for all \( t \in \mathbb{R} \), so \( \sigma \) vanishes identically on the line \( D := \{ t\zeta : t \in \mathbb{R} \} \).

Let \( \{ \xi_j \}_{j=1}^\infty \subset \mathbb{Z}^d \) be a sequence of lattice points with \( \operatorname{dist}(\xi_j, D) \lesssim 1 \) and \( |\xi_j| \to \infty \). For each \( j \), let \( \zeta_j \in D \) denote the orthogonal projection of \( \xi_j \) onto \( D \), and write
\[
\zeta_j^\perp := \frac{\xi_j - \zeta_j}{|\xi_j - \zeta_j|}.
\]
Define \( F_j(t) := \sigma(\zeta_j + t \zeta_j^\perp) \), which is a degree-\( m \) polynomial in \( t \). Since \( \sigma(\zeta_j) = 0 \), we have:
\[
|\sigma(\xi_j)| = |F_j(|\xi_j - \zeta_j|)| \leq \sum_{k = 1}^m \frac{|F_j^{(k)}(0)|}{k!} |\xi_j - \zeta_j|^k.
\]
Each derivative \( F_j^{(k)}(0) \) is a sum of terms involving monomials \( \zeta_j^\alpha \), where \( \zeta_j \in D \) and \( |\zeta_j| \sim |\xi_j| \). Therefore:
\[
|F_j^{(k)}(0)| \lesssim |\zeta_j|^{m - k} \lesssim |\xi_j|^{m - k},
\]
and since \( |\xi_j - \zeta_j| \lesssim 1 \), we obtain:
\[
|\sigma(\xi_j)| \lesssim \sum_{k = 1}^m |\xi_j|^{m - k} \lesssim |\xi_j|^{m - 1}.
\]
Finally, for large \( j \), we have:
\[
|P(\xi_j) - \lambda| \leq |\sigma(\xi_j)| + |r(\xi_j)| + |\lambda| \lesssim |\xi_j|^{m - 1}.
\]
Thus, for any \( \delta > 0 \), we have:
\[
|P(\xi_j) - \lambda| \leq |\xi_j|^{m - 1 + \delta}
\]
for all sufficiently large \( j \), showing that the Diophantine inequality has infinitely many integer solutions, a contradiction. This completes the proof.
\end{proof}

We have the following bound on the lattice point count in the case of real coefficients.
\begin{lemma}[Asymptotic Estimate for the Lattice Point Count] \label{lem:lattice-count}
Let $P(\xi)$ be the symbol of an elliptic differential operator of order $m \geq 1$ with 
{\bf real} coefficients, and fix $\delta > 0$. Define
\[
\mathcal{N}_\lambda := \#\left\{ \xi \in \mathbb{Z}^d : |P(\xi) - \lambda| \le |\xi|^{m-1+\delta} \right\}.
\]
Then, as $|\lambda| \to \infty$,
\[
\mathcal{N}_\lambda \lesssim |\lambda|^{\frac{d-1+\delta}{m}},
\]
where the implitict constant depends on $P$, $d$, and $\delta$.
\end{lemma}

\begin{proof}
Fix large $|\lambda|$ and observe that the set
\[
R(\lambda) := \left\{ \xi \in \mathbb{R}^d : |P(\xi) - \lambda| \le |\xi|^{m-1+\delta} \right\}
\]
forms a thickened neighborhood around the level set $\{ P(\xi) = \lambda \}$.

Denote by $P_m$ the principal symbol of $H$. By definition we have  $P(\xi) = P_m(\xi) + O(|\xi|^{m-1})$, as $|\xi| \to \infty$, so we may rewrite the inequality as
\begin{equation}\label{dioph1}
\Big|P_m(\frac{\xi}{|\xi|}) +O(\frac{1}{|\xi|})- \frac{\lambda}{|\xi|^m}\Big| \le |\xi|^{-1+\delta}
\end{equation}

Because of ellipticity we know that $P_m$ cannot vanish on $S^{d-1}$ and therefore does not change sign. Without loss of generality, let us assume that $P_m>0$ on $S^{d-1}$. Moreover, assume $\lambda>0$. We wish to estimate the number of solutions to the above diophantine inequality.

Note that \eqref{dioph1} implies 
\[
\Big|P_m(\theta)- \frac{\lambda}{|\xi|^m}\Big| \lesssim |\xi|^{-1+\delta}
\]
where $\theta:=\xi/|\xi|$ denotes the projection of $\xi$ onto $S^{d-1}$.

We rewrite the inequality as
\[
\left( \frac{\lambda}{P_m(\theta)+C|\xi|^{-1+\delta}}\right)^{1/m} \leq |\xi| \leq \left(\frac{\lambda}{P_m(\theta)-C|\xi|^{-1+\delta}}\right)^{1/m}
\]
for some $C>0$. Because for any $\theta\in S^{d-1}$ we have $P_m(\theta)\in[c_1,c_2]$ for some constants $0<c_1<c_2$, it follows that $|\xi|\asymp \lambda^{1/m}$.

Hence, we calculate the width of the region described by this inequality as
\[w=y^{1/m}-x^{1/m}\]
where \[x:=\frac{\lambda}{P_m(\theta)+C|\xi|^{-1+\delta}} \;\text{and}\; y:=\frac{\lambda}{P_m(\theta)-C|\xi|^{-1+\delta}}.\]

We obtain 
\[
w\asymp c^{1/m-1}(y-x)\asymp \lambda^{1/m-1}\lambda^{1+1/m\cdot (-1+\delta)}=\lambda^{\delta/m},
\]
where $c\in(x,y)$ and we used $x\asymp y\asymp \lambda$.

Hence, the width grows with $\lambda$ and the region described by the above diophantine inequality is a distorted spherical shell of width $O(\lambda^{\delta/m})$ and central radius $\asymp \lambda^{1/m}$. The number of lattice points in this region is thus estimated by the number of lattice points in an annulus of central radius $R=\lambda^{1/m}$ and width $w=\lambda^{\delta/m}$ which is $O(R^{d-1}w)=O(\lambda^{(d-1+\delta)/m})$.
\end{proof}

In the case of complex coefficients it is easy to find operators which admit much better bounds.
\begin{lemma}\label{lem:complex-count}
Let $H=-\partial_{x}^2-i\partial_{y}^2$. In this case, we have  
\begin{equation}
\mathcal{N}_\lambda\lesssim
\begin{cases} 
(\re\lambda)^{1/2+\delta}(\im\lambda)^{-1/2}\quad \text{if}\; (\re\lambda)^{1/2}\leq\im\lambda\leq\re\lambda,\\
\\
(\re\lambda)^{1/4+\delta} \quad\text{if}\; \im\lambda\leq(\re\lambda)^{1/2}.
\end{cases}
\end{equation}
if $0<\im\lambda\leq\re\lambda\to+\infty$. Analogous bounds hold if one exchanges the roles of $\re\lambda$ and $\im\lambda$.
\end{lemma}
\begin{proof}
The symbol of the operator is $P(\xi)=\xi_1^2+i\xi_2^2$. Take $\lambda=a+ib$, $a,b>0$. The inequality $$|P(\xi)-\lambda|\leq |\xi|^{1+\delta}$$ implies
$$|\xi_1^2-a|\leq |\xi|^{1+\delta}, \; \text{and}\; |\xi_2^2-b|\leq |\xi|^{1+\delta}.$$

Assume $|\xi_2|\leq|\xi_1|$. Then
$$|\xi_1^2-a|\lesssim |\xi_1|^{1+\delta}$$ and the solutions $\xi_1\in\Z$ are located in an interval centered on $a^{1/2}$ of width $O(a^{\delta/2})$.
Now, let us estimate the number of solutions $\xi_2\in\Z$ to the second inequality
$$|\xi_2^2-b|\leq |\xi|^{1+\delta}\lesssim a^{(1+\delta)/2}$$
which is bounded by $O\left(\frac{a^{(1+\delta)/2}}{(b+a^{(1+\delta)/2})^{1/2}}\right)$. So the overall number of solutions $\xi=(\xi_1,\xi_2)$ is bounded by $$O\left(\frac{a^{1/2+\delta}}{(b+a^{(1+\delta)/2})^{1/2}}\right).$$

If $|\xi_1|\leq|\xi_2|$, we follow the same argument and obtain the bound $$O\left(\frac{b^{1/2+\delta}}{(a+b^{(1+\delta)/2})^{1/2}}\right).$$

To conclude, let us suppose $0<b\leq a$. If $0<b\leq a^{1/2}$, we get $O(a^{1/4+\delta})$. On the other hand, if $a^{1/2}\leq b\leq a$, then we get $O((a/b)^{1/2}a^{\delta})$. 
\end{proof}


\appendix


\section{Representation of Integers as Sums of Squares} \label{sec:rdn}

We denote by \( r_d(n) \) the number of ways a positive integer \( n \) can be written as a sum of \( d \) integer squares, i.e.,
\[
r_d(n) := \#\left\{ \xi \in \mathbb{Z}^d : |\xi|^2 = n \right\},
\]
which is the number of lattice points on the \( (d-1) \)-sphere of radius \( \sqrt{n} \) in \( \mathbb{R}^d \).

\subsection{Case \( d = 2 \)}

For two dimensions, this count is known to grow very slowly:
\[
r_2(n) \lesssim e^{c \log n / \log \log n}, \quad \text{as } n \to \infty,
\]
for some constant \( c > 0 \). This reflects the arithmetic sparsity of representations of integers as sums of two squares.

\subsection{Case \( 3 \le d \le 8 \)}

In dimensions \( 3 \le d \le 8 \), Hardy \cite{Ha1920} gave an exact formula:
\[
r_d(n) = \frac{\pi^{d/2}}{\Gamma(d/2)} n^{d/2 - 1} S_d(n),
\]
where \( S_d(n) \) is the \emph{singular series}, given by:
\[
S_d(n) = \sum_{k = 1}^\infty k^{-1/2} \sum_{h = 1}^{2k} \eta(h,k)^d e^{-\pi i hn / k},
\]
with
\[
\eta(h,k) :=
\begin{cases}
\frac{1}{2}k^{-1/2} \sum_{j = 1}^{2k} e^{\pi i h j^2 / k} & \text{if } \gcd(h, k) = 1, \\
0 & \text{otherwise}.
\end{cases}
\]

The singular series captures the arithmetic structure of \( n \) and may exhibit mild fluctuations. For \( d \ge 5 \), however, it is uniformly bounded.

\subsection{Case \( d > 8 \)}

For higher dimensions, Hardy and Ramanujan \cite{Ra18,Ha1920} proved the following asymptotic:
\[
r_d(n) = \frac{\pi^{d/2}}{\Gamma(d/2)} n^{d/2 - 1} S_d(n) + O(n^{d/4}), \quad \text{as } n \to \infty.
\]
Here the main term reflects the surface area of the sphere, while the error term accounts for finer arithmetic effects. The error is negligible for large \( d \).

\subsection{Growth of the Singular Series}

From analytic number theory (see e.g., Bateman \cite{Bt51}), the singular series satisfies the following bounds:
\[
S_d(n) \lesssim
\begin{cases}
\log n \cdot \log \log n & \text{if } d = 3, \\
\log \log n & \text{if } d = 4, \\
1 & \text{if } d \ge 5.
\end{cases}
\]
Hence, for \( d \ge 3 \), we obtain:
\[
r_d(n) \lesssim n^{d/2 - 1} \times
\begin{cases}
\log n \cdot \log \log n & \text{if } d = 3, \\
\log \log n & \text{if } d = 4, \\
1 & \text{if } d \ge 5.
\end{cases}
\]

These estimates ensure that the series
\[
\sum_{n = 1}^\infty r_d(n) n^{-s}
\]
converges for \( \operatorname{Re}(s) > d/2 \), and justify the decay rates needed in the lattice point tail estimates used throughout the paper.

\section{Shrinking Balls} \label{sec:shrink}

In this section, we refine Lemma~\ref{lem:fourier-bound} in the special case \( H = -\Delta \), in order to make explicit the dependence of the bound on the radius \( r > 0 \). This allows us to analyze the effect of allowing \( r = r(\lambda) \) to decay slowly with \( \lambda \), which is essential in rough or unbounded domains.

Recall that the cutoff is defined by \( \chi_r(x) := \chi(x/r) \), so for any multi-index \( \gamma \) we have:
\[
D^\gamma \chi_r(x) = r^{-|\gamma|} (D^\gamma \chi)(x/r).
\]

As shown in the proof of Lemma~\ref{lem:fourier-bound}, the quantity \( (-\Delta - \lambda)^\alpha \Psi \) is a sum of terms of the form \( h \cdot D^\gamma \psi_\lambda \), where \( h \in C^\infty_c(B(x_0, r)) \) depends on derivatives of \( \chi_r \) and \( V \), and \( |\gamma| \leq 2\alpha - \deg(h) \).

Thus, each such term satisfies
\[
|h(x)| \lesssim r^{-k}, \quad \text{for some } k \leq 2\alpha - |\gamma|.
\]

We then obtain the following refined estimate for the Fourier coefficients:
\[
\big||\xi|^2 - \lambda \big|^\alpha |\hat{\Psi}(\xi)| \lesssim_{\chi,\alpha} \sum_{k = 0}^\alpha r^{-2\alpha + k} |\xi|^k \|\psi_\lambda\|_{L^2(B(x_0, r))}.
\]

\subsection{Dependence on \( r \) in Part II of the Summation}

We now revisit Part II of the Fourier sum in the proof of Theorem~\ref{thm:local-bound}, using the refined bound above. For the region where
\[
|P(\xi) - \lambda| = ||\xi|^2 - \lambda| \geq |\xi|^{1+\delta},
\]
we estimate:
\begin{align*}
\sum_{||\xi|^2 - \lambda| \geq |\xi|^{1+\delta}} |\hat{\Psi}(\xi)|
&\lesssim \sum_{k = 0}^\alpha r^{-2\alpha + k} \sum_{\xi \neq 0} \frac{|\xi|^k}{|\xi|^{(1+\delta)\alpha}} \|\psi_\lambda\|_{L^2(B(x_0, r))} \\
&= \|\psi_\lambda\|_{L^2(B(x_0, r))} \sum_{k = 0}^\alpha r^{-2\alpha + k} \sum_{\xi \neq 0} |\xi|^{- \alpha \delta + (k - \alpha)}.
\end{align*}

The inner sum converges if \( \alpha\delta > d + k - \alpha \), which is ensured by taking \( \alpha = \alpha(\delta) > d / \delta \). Therefore, we obtain the uniform bound:
\[
\sum_{||\xi|^2 - \lambda| \geq |\xi|^{1+\delta}} |\hat{\Psi}(\xi)| \lesssim_{d,\delta} r^{-2\alpha} \|\psi_\lambda\|_{L^2(B(x_0, r))}.
\]

\subsection{Final Estimate with Explicit Radius Dependence}

Combining this with the bound from Part I, we obtain:
\[
|\psi_\lambda(x)| \lesssim \left( r^{-2\alpha(\delta)} + \lambda^{(d - 1 + \delta)/4} \right) \|\psi_\lambda\|_{L^2(B(x_0, r))}.
\]

Hence, if we allow the radius to decay slowly with \( \lambda \), say
\[
r^{-1} = \lambda^{o(1)} \quad \text{as } \lambda \to \infty,
\]
then the contribution from the shrinking-ball term remains negligible compared to the main term \( O_\varepsilon(\lambda^{(d - 1)/4 + \varepsilon}) \) when choosing \( \delta = 4\varepsilon \).

This justifies the use of shrinking interior balls \( B(x_0, r(\lambda)) \) in Corollary~\ref{cor:interior-bound}, even for large spectral parameters \( \lambda \).


\end{document}